

\frenchspacing

\parindent15pt

\abovedisplayskip4pt plus2pt
\belowdisplayskip4pt plus2pt 
\abovedisplayshortskip2pt plus2pt 
\belowdisplayshortskip2pt plus2pt  

\font\twbf=cmbx10 at12pt
 at12pt
 at12pt

\font\sc=cmcsc10

\font\ninerm=cmr9 
\font\nineit=cmti9 
\font\ninesy=cmsy9 
\font\ninei=cmmi9 
\font\ninebf=cmbx9 

\font\sevenrm=cmr7  
 
\font\seveni=cmmi7  
\font\sevensy=cmsy7 

\font\fivenrm=cmr5  
\font\fiveni=cmmi5  
\font\fivensy=cmsy5 

\def\nine{%
\textfont0=\ninerm \scriptfont0=\sevenrm \scriptscriptfont0=\fivenrm
\textfont1=\ninei \scriptfont1=\seveni \scriptscriptfont1=\fiveni
\textfont2=\ninesy \scriptfont2=\sevensy \scriptscriptfont2=\fivensy
\textfont3=\tenex \scriptfont3=\tenex \scriptscriptfont3=\tenex
\def\rm{\fam0\ninerm}%
\textfont\itfam=\nineit    
\def\it{\fam\itfam\nineit}%
\textfont\bffam=\ninebf 
\def\bf{\fam\bffam\ninebf}%
\normalbaselineskip=11pt
\setbox\strutbox=\hbox{\vrule height8pt depth3pt width0pt}%
\normalbaselines\rm}

\hsize30cc
\vsize44cc
\nopagenumbers

\def\luz#1{\luzno#1?}
\def\luzno#1{\ifx#1?\let\next=\relax\yyy
\else \let\next=\luzno#1\xxx\fi\next}
\def\sp#1{\def\xxx{\kern1.7pt}\def\yyy{\kern-1.7pt}\luz{#1}}
\def\spa#1{\def\xxx{\kern1pt}\def\yyy{\kern-1pt}\luz{#1}}

\newcount\beg
\newbox\aabox
\newbox\atbox
\newbox\fpbox
\def\abbrevauthors#1{\setbox\aabox=\hbox{\sevenrm\uppercase{#1}}}
\def\abbrevtitle#1{\setbox\atbox=\hbox{\sevenrm\uppercase{#1}}}
\long\def\pag{\beg=\pageno
\def\leftheadline{\noindent\rlap{\nine\folio}\hfil\copy\aabox\hfil}
\def\rightheadline{\noindent\hfill\copy\atbox\hfill\llap{\nine\folio}}
\def\phead{\setbox\fpbox=\hbox{\sevenrm 
************************************************}%
\noindent\vbox{\sevenrm\baselineskip9pt\hsize\wd\fpbox%
\centerline{***********************************************}

\centerline{BANACH CENTER PUBLICATIONS, VOLUME **}

\centerline{INSTITUTE OF MATHEMATICS}

\centerline{POLISH ACADEMY OF SCIENCES}

\centerline{WARSZAWA 19**}}\hfill}
\footline{\ifnum\beg=\pageno \hfill\nine[\folio]\hfill\fi}
\headline{\ifnum\beg=\pageno\phead
\else
\ifodd\pageno\rightheadline \else \leftheadline \fi 
\fi}}

\newbox\tbox
\newbox\aubox
\newbox\adbox
\newbox\mathbox

\def\title#1{\setbox\tbox=\hbox{\let\\=\cr 
\baselineskip14pt\vbox{\twbf\tabskip 0pt plus15cc
\halign to\hsize{\hfil\ignorespaces \uppercase{##}\hfil\cr#1\cr}}}}

\newbox\abbox
\setbox\abbox=\vbox{\vglue18pt}

\def\author#1{\setbox\aubox=\hbox{\let\\=\cr 
\nine\baselineskip12pt\vbox{\tabskip 0pt plus15cc
\halign to\hsize{\hfil\ignorespaces \uppercase{\spa{##}}\hfil\cr#1\cr}}}%
\global\setbox\abbox=\vbox{\unvbox\abbox\box\aubox\vskip8pt}}

\def\address#1{\setbox\adbox=\hbox{\let\\=\cr 
\nine\baselineskip12pt\vbox{\it\tabskip 0pt plus15cc
\halign to\hsize{\hfil\ignorespaces {##}\hfil\cr#1\cr}}}%
\global\setbox\abbox=\vbox{\unvbox\abbox\box\adbox\vskip16pt}}

\def\mathclass#1{\setbox\mathbox=\hbox{\footnote{}{1991 {\it Mathematics Subject 
Classification}\/: #1}}}

\long\def\maketitlebcp{\pag\unhbox\mathbox
\footnote{}{The paper is in final form and no version 
of it will be published elsewhere.} 
\vglue7cc
\box\tbox
\box\abbox
\vskip8pt}

\long\def\abstract#1{{\nine{\bf Abstract.} 
#1

}}

\def\section#1{\vskip-\lastskip\vskip12pt plus2pt minus2pt
{\bf #1}}

\long\def\th#1#2#3{\vskip-\lastskip\vskip4pt plus2pt
{\sc #1} #2\hskip-\lastskip\ {\it #3}\vskip-\lastskip\vskip4pt plus2pt}

\long\def\defin#1#2{\vskip-\lastskip\vskip4pt plus2pt
{\sc #1} #2 \vskip-\lastskip\vskip4pt plus2pt}

\long\def\remar#1#2{\vskip-\lastskip\vskip4pt plus2pt
\sp{#1} #2\vskip-\lastskip\vskip4pt plus2pt}

\def\Proof{\vskip-\lastskip\vskip4pt plus2pt 
\sp{Proo{f.}\ }\ignorespaces}

\def\endproof{\nobreak\kern5pt\nobreak\vrule height4pt width4pt depth0pt
\vskip4pt plus2pt}

\newbox\refbox
\newdimen\refwidth
\long\def\references#1#2{{\nine
\setbox\refbox=\hbox{\nine[#1]}\refwidth\wd\refbox\advance\refwidth by 12pt%
\def\textindent##1{\indent\llap{##1\hskip12pt}\ignorespaces}
\vskip24pt plus4pt minus4pt
\centerline{\bf References}
\vskip12pt plus2pt minus2pt
\parindent=\refwidth
#2

}}

\def\footnoterule{\kern -3pt \hrule width 4cc \kern 2.6pt}

\catcode`@=11
\def\vfootnote#1%
{\insert\footins\bgroup\nine\interlinepenalty\interfootnotelinepenalty%
\splittopskip\ht\strutbox\splitmaxdepth\dp\strutbox\floatingpenalty\@MM%
\leftskip\z@skip\rightskip\z@skip\spaceskip\z@skip\xspaceskip\z@skip%
\textindent{#1}\footstrut\futurelet\next\fo@t}
\catcode`@=12


\def\ba{\begin{array}}
\def\ea{\end{array}}
\def\lra{\longrightarrow}
\def\ra{\rightarrow}
\def\bk{I \kern-.25em k} 
\def\lin{\hbox{lin}\,}
\def\der{\hbox{der}\,}
\def\derk{\der_{\bk}\,}

\def\im{\hbox{Im}\,}
\def\ker{\hbox{ker}\,}
\def\ann{\hbox{ann}\,}
\def\bot{\,_b\!\otimes}
\def\nd{\hbox{End}\,}
\def\hom{\hbox{hom}\,}
\def\alg{\hbox{alg}\,}

\def\id{\hbox{id}}
\def\idA{\id_A}
\def\idV{\id_V}
\def\ot{\otimes}
\def\c{\circ}

\mathclass{Primary 16D20; Secondary 16W25; 16U80.}

\abbrevauthors{A. Borowiec, V. Kharchenko and Z. Oziewicz}
\abbrevtitle{Calculi with universal bimodule}

\title{First order calculi with values\\
in  right--universal bimodules}

\author{Andrzej\ Borowiec}
\address{Institute of Theoretical Physics, University of Wroc{\l}aw\\
Pl. Maxa Borna 9, 50-204 Wroc{\l}aw, Poland\\
E-mail: borow@ift.uni.wroc.pl}
\author{Vladislav\ Kharchenko}
\address{Centro de Investigaciones Teoricas, FESC, UNAM,\\ 
Apartado Postal 95,
C.P. 54700 Cuautitl\'an Izcalli, Estado de M\'exico.\\
E-mail: vlad@servidor.unam.mx}
\author{Zbigniew\ Oziewicz}
\address{Facultad de Estudios Superiores,
Cuautitl\'an,\\Universidad Nacional Autonoma de M\'exico,\\Apartado
Postal 25, C.P. 54700 Cuautitl\'an Izcalli, Estado de M\'exico.\\
E-mail: oziewicz@servidor.unam.mx}

\maketitlebcp

\footnote{}{Research partially supported by KBN grant 2 P302 023 07.}
\footnote{}{The second author wishes to thank CONASyT--Mexico for its
support, Catedra Patrimonial 

940411 and also the Russian Foundation
for Fundamental Research, grant 95-01-01356.}
\footnote{}{Z. Oziewicz is a member of Sistema Nacional de Investigadores, 
M\'exico.}

\abstract{The purpose of this note is to show how calculi on unital
associative algebra with universal right bimodule generalize
previously studied constructions by Pusz and Woronowicz [1989] and by
Wess and Zumino [1990] and that in this language results are in a
natural context, are easier to describe and handle.
As a by--product we obtained intrinsic, coordinate--free and 
basis--independent generalization of the first order 
noncommutative differential calculi with partial derivatives.}

\section{1. Introduction.}
In this note $\bk$ denotes some fixed unital and commutative ring. 
Algebras are unital associative $\bk$-algebras and homomorpisms
are assumed to be unital. The tensor product $\ot$ 
means $\ot_{\bk}$ and linear map
means $\bk$-linear. All objects considered here  are first of all 
$\bk$-modules. 
All maps are assumed to be $\bk$-linear maps. Most results of
Section 3 will be set under additional assumption that
$\bk$-modules $V$ and $W$ are projective.
\vskip4pt plus2pt

{\bf 1.1.} {\it First order differential calculus as a derivation
of an algebra.}
Let $A$ be an algebra with a multiplication
$m\in\lin(A^{\otimes 2},A)$, $M$ be an
$(A, A)$-bimodule ($A$-bimodule in short).

It  appears that differential calculi investigated recently by many
authors in the context of quantum groups and noncommutative 
geometry are nothing but dervations of an algebra with values in
a bimodule.
Recall that a $\bk$-derivation $d$ of $A$ to $M,$
$d\in\der_{\bk}(A, M)$, is a $\bk$-linear mapping from
$A$ into $M$ such that the Leibniz rule
$$
d(x y) = dx.y + x.dy 
$$
holds true for each $x,\,y\in A$ where, dot "." denotes both 
(left and right ) multiplications by elements from $A$.
\defin{Definition}{1. 
The triple $\{A, M, d\}$ is said to be {\it first order calculus} or 
{\it first order differential} on an algebra $A$ with values in 
an $A$--bimodule $M$ or shortly {\it $M$-valued calculus} on $A$.}
Each $\bk$-derivation vanishes on scalars from $\bk$.
\vskip4pt plus2pt
{\bf 1.2.} {\it Functorial properties.}
The following functorial properties of derivations are well known [2]:
\th{Proposition}{1.2.1. (Bourbaki 1989, \S 10.7, Prop. 9).}{
Let $A,\ B$ be two algebras, $M$ an $A$-bimodule and $N$ a 
$B$-bimodule;
let $\phi :A\lra B$ be an algebra homomorphism and $\Phi : M\lra N$ 
an $A$-homomorphism of $A$-bimodules (relative to $\phi$). Then
\item{(i)} For every $d_B\in \derk (B, N)$, $d_B\c\phi\in \derk 
(A, \phi^{\ast}(N))$.
\item{(ii)} For every $d_A\in \derk (A, M)$, $\Phi\c d_A\in \derk 
(A, \phi^{\ast}(N))$.
where, $\phi^{\ast}(N)=N$ with an $A$-bimodule structure  induced by $\phi$.
}
We are interested in the case when the diagram
$$
\matrix{A &{\buildrel d_A\over\longrightarrow} & M\cr
\phi \downarrow &   & \downarrow \Phi\cr
B &{\buildrel d_B\over\longrightarrow} & N 
}\leqno(1)$$
commute. The condition of commutativity 
$$\Phi\c d_A=d_B\c \phi \leqno (2)$$
does not in general determine uniquely the derivations $d_A, \,d_B$. 
If $\phi$ is surjective then for a given $(\phi, \, \Phi, \, d_A)$
there exists at most one $d_B$ rendering commutative the diagram (1).
The existence problem will follow from an additional assumption
so called {\it consistency conditions}. 

In fact we will be more specific.
Assume that $\phi, \Phi$ are surjective homomorphisms and the diagram (1)
commute. One has $B\cong A/\ker\phi$ and $N\cong M/\ker\Phi$ {\it i.e.}
$d_B$ (if exists) is a {\it factor calculus} on factor algebra 
with values in factor bimodule.
The differential $d_A$ is called a {\it cover differential}.
In the next Section we shall study factor calculi and the corresponding
consistency conditions with more details.
\th{Lemma}{1.2.2 ([2], p. 560) }{
Let $V$ be a generating system for the algebra $A$.
The diagram (1) is commutative if and only if it is commutative
on elements from $V$.
}

{\bf 1.3.} {\it Universal problem for calculi.}
It is also  very well known that for any algebra $A$ there exists an universal
derivation $\delta_A\in\derk (A, \ker m)$ with values in a kernel of 
multiplication map $m$ and $\delta_A$ is defined by the formulae
$$
\delta_A= \idA\ot 1_A - 1_A\ot\idA. \leqno(3)
$$
with $1_A$ being a unit of $A$.
It has the following universal property [2, 13]: for every 
$A$-bimodule $M$ and every derivation $d\in\derk (A, M)$, there exists
one and only one $A$-bimodule homomorphism $\Phi_d :\ker m\ra M$
such that $d=\Phi_d\c\delta_A$. In other words an universal differential
is a cover differential for any derivation $d$ of the algebra $A$ 
(one has $\phi_d=\idA$ in this case). The existence is proved by setting
$$
\Phi_d(\sum_i x_i\ot y_i) =\sum_i dx_i .y_i\equiv -\sum_i x_i .dy_i
$$
where, $\sum_i x_i\ot y_i\in\ker m$. Thus we have  a canonical
$\bk$-module isomorphism
$$
\hom_{(A, A)}(\ker m, M)\ {\buildrel \delta_A\over\longleftrightarrow}\ 
\derk (A, M) . 
$$ 
\remar{Remark\ {1.}\ }{
Note that $\ker m$ is generated (as a left or right $A$-module)
by the image of $\delta_A$.
Therefore if $M$ is also generated by differentials {\it i.e.}
$M=dA.A=A.dA$ then $M$ is a factor bimodule of $\ker m$
by some subbimodule. This condition was included by Woronowicz into
the definition of calculi [18]. 
In this case $\Phi_d$ is surjective.}

{\bf 1.4\ }{\it Calculi with partial derivatives.}
Assume that $\bk$ is a field.
The algebra $A$ possesses a basis $\{ 1_A, \,e_s\}_{s\in S}$ 
(consisting finite or infinite number of elements) as a linear
$\bk$-space.
One can easily observe (see {\it e.g.} [5]) that the set 
$\{e_s\ot 1-1\ot e_s\}_{s\in S}$
is a set of free generators for the right (or left) $A$-module $\ker m$.
The condition that bimodule $M$ of one-forms admits free generators
was assumed in the celebrated paper by Pusz and Woronowicz [14].
It was also a guiding line (although not assumed explicitly)
in the paper by Wess and Zumino [17].

Let $d\in \derk (A, M)$ and assume that the bimodule $M$ as a
right $A$-module is free. Let elements $\mit\{ \xi^\gamma\}_{
\gamma\in\Gamma}$ be free generators of $M$ as a right module. 
One has in this case
$$\mit
da=\sum_{\gamma\in\Gamma}\xi^\gamma D_\gamma a 
$$
where, the {\it partial derivatives} $\mit D_\gamma\in \hom_{\bk}(A, A)$ 
are uniquely defined. We call such calculus a {\it calculus with partial 
derivatives}. In particular, if $\mit\{ x^\gamma\}_{\gamma\in\Gamma}$
is a set of generators of the algebra $A$ such that 
$\mit dx^\gamma=\xi^\gamma$
then the corresponding calculus is called {\it coordinate calculus}.
Recently there is a great interest in non commutative differential
calculi mainly with connection to quantum groups and covariant 
calculi on quantum hyperplains [9--12] and [14--18].
Many examples of coordinate calculi and calculi with partial
derivatives have been studied in the literature (see {\it e.g.} 
[1, 5, 6, 9, 15]).
More systematic approach and general formalism for coordinate calculi 
and calculi with partial derivatives over a field have been proposed
[3, 4].

Let $W$ denotes a free $\bk$-module with  basis 
$\mit\{ \xi^\gamma\}_{\gamma\in\Gamma}$. $W\ot A$ has a canonical
right $A$-module structure and 
with this structure it is isomorphic (but not canonically) to 
$M$ when is considered as a right $A$-module.

Let $V$ be a generating $\bk$-module for the algebra $A$ .
The aim of the present paper is to describe  differential calculi on
bimodules of the type $W\ot A$ where, $V,\, W$ are some $\bk$-modules.
We do not assume in general that $V,\, W$ are free modules.
In this way one obtains coordinate--free and 
basis--independent generalization of calculi with partial derivatives
as well as a generalization to the case when $\bk$ is a commutative
ring. This will be done in Section 3.

Our approach to non--commutative first--order differential calculi
is traditionally a covariant one: it is done by means of differential
one--forms. An extension to forms of higher degree 
(quantum de Rham complexes) can be found {\it e.g.}
in [4, 10--12, 14, 17, 18]. A contravariant ({\it i.e.} by means of 
vector fields) approach will be proposed in [7].

\section{2. Factor calculus on factor algebra with factor bimodule.}
Let now $A, B$ be two algebras, $M$ an $A$-bimodule and $N$
a $B$-bimodule. Notice that any homomorphism $\phi :A\ra B$ induces
an $A$-bimodule structure on $N$ (change of scalars) by
$$
x.n=\phi (x).n ,\ \ n.x=n.\phi (x) \ \ \ .$$
We  denote this new $A$-bimodule structure on $N$ by $\phi^\ast (N)$
($\phi^\ast (N)=N$  as sets).

In the rest of this section $A, \,B$ are $\bk$ algebras, $M$ 
is an $A$-bimodule
and $\phi :A\ra B$ is an surjective homomorphism.
Recall that annihilator of a subbimodule $X\subseteq M$ in $A$ is a
two sided ideal $\ann X$ defined by
$$
\ann X = \{a\in A: a.X=X.a=0\} \ \ \leqno(4)
$$
\vskip4pt plus2pt
{\bf 2.1.} {\it Factor bimodule.}
In what follows we shall need the obvious criterion for a factor
module $M/L$ to be an $A$-bimodule.
\th{Lemma}{2.1.1. } {Let $L$ be a $\bk$-submodule of the $A$-bimodule $M$. 
The following are equivalent:
\item{(i)} $L$ is an $A$-subbimodule of $M$.
\item{(ii)} $M/L$ with the induced multiplication is an $A$-bimodule.
\item{(iii)} $A.L.A\subseteq L$.
\item{(iv)} $A.L + L.A\subseteq L$.
}
One sees that for given $\bk$-submodule $L$, $\lin (A.L.A)$ is the smallest
$A$-subbimodule of $M$ containing $L$.

\th{Proposition}{2.1.2. (Pierce 1982) }
{Let $I=\ker \phi$. Then $M=\phi^\ast N$ for some
$B$-module $N$ if and only if 
$I\subseteq \ann M \ \ $.
}
\Proof For a case of right modules the proof can be found in ref. [13] p.23 .
\vskip4pt plus2pt
In other words we get a condition for the $A$-bimodule $M$
to become a $B\cong A/I$-module.
Below we will make no distinction between $A/I$-bimodules and $A$-bimodules
$M$ if $I\subseteq \ann M$. The ideal $\ann M$ is a maximal ideal
possessing this property. A faithful bimodule ($\ann =0$) cannot be 
a bimodule over proper factor algebra in a way as describe above.

\th{Proposition}{2.1.3 }{ Let $L$ be a subbimodule of $M$ and 
$I < A$ an ideal of $A$.
Then the structure of $A$-bimodule on $M$ induces a structure of
$A/I$-bimodule if and only if 
$$I.M+M.I\subseteq L \ \ .\leqno (5)$$
}
\Proof From the previous Proposition we have to have $I\subseteq \ann M/L$,
which is equivalent to $I.(M/L)=(M/L).I=0$. It gives the proof.
\vskip4pt plus2pt

One sees that $\lin (I.M+M.I)$ is a subbimodule of $M$. It is the smallest
subbimodule satisfying the statement above. We shall call the condition
(5) {\it compatibility condition} between an ideal $I$ and a subbimodule $L$.
Therefore, $\lin (I.M+M.I)$ is the smallest $I$-compatible subbimodule of $M$.
If $I$ and $L$ are compatible then the factor module $M/L$ becomes a
module over factor algebra $A/I$. In particular, if $I\subseteq \ann M$
{\it i.e.} $M.I=I.M=0$, then any subbimodule $L$ is $I$-compatible.
Observe that for a faithful bimodule $M$, $M/\lin (I.M+M.I)$ is in
some sense a maximal factor bimodule over factor algebra $A/I$.
From the other side if a subbimodule $L$ is fixed then there exists
the largest $L$-compatible ideal $I(L)$ 
(the sum of all ideals satisfying (5)). 
The algebra $A/I(L)$ is in some sense a
minimal algebra compatible with the subbimodule $L$.
\vskip4pt plus2pt
{\bf 2.2.} {\it Factor calculus on factor algebra with factor bimodules.}
Combining Proposition 1.2.1 with the results above we have

\th{Theorem}{2.2.1. (Factor calculi)} {Let $ I < A$ be an ideal
in $A$, $L$ be a subbimodule of an $A$-bimodule $M$ and 
$d\in\derk (A, M)$. In order to exist a unique factor calculus
$\tilde d\in\derk (A/I, \,M/L)$ such that
$$\Pi\c d=\tilde d\c\pi \leqno(6)$$
where, $\pi :A\ra A/I$ and $\Pi :M\ra M/L$ are canonical surjections,
it is sufficient and necessary that
$$dI + I.M + M.I \subseteq L \ \ . \leqno (7)$$
}
\Proof We know that $M.I+I.M\subseteq L$ guarantees existence of
an $A/I$ bimodule structure on $M/L$.
One also knows that $\tilde d$ if exists is unique. Moreover
$\tilde d(x+I)=dx+L$ is well defined if and 
only if  $dI\subseteq L$. This leads to commutativity condition (6).
To complete the proof one has to apply the Proposition 1.2.1 .
\vskip4pt plus2pt
\remar{Remark\ {2.}}{Notice that in the case $I=0$ {\it i.e.}
$\pi=\idA$ there is no any condition for existence of factor calculi
on a factor bimodule $M/L$.
}

When $L$ is a $\bk$-module then the condition (7) must 
be replaced by (c.f. Lemma 2.1.1)
$$
dI+M.I+I.M + A.L.A\subseteq L \ .\leqno(8)
$$

The conditions of the type (7, 8) for existence a factor derivation are 
usually in the literature called {\it consistency conditions} for 
a factor calculi.
They include a compatibility condition between an ideal $I$ and
$A$-subbimodule ($\bk$-submodule) $L$.
Notice that in this case the cover differential $d$ is not in 
general, uniquely determined by the diagram (1).

Again one can consider some special situations.
Firstly, remark that if $I\subseteq \ker d$ then (5) is the only condition 
ensuring the existence of factor calculus. Secondly, for a given $I$, 
$\lin (dI+M.I+I.M)$ is the smallest subbimodule admitting a factor
calculus. Finally if $L$ is fixed then there exists a maximal
ideal $I(d, L)$ and the corresponding optimal algebra with a 
factor calculus. Notion of optimal algebra has been introduced
in [3] (c.f. [4--6]).

\section{3. Calculi with right--universal bimodules}

\vskip4pt plus2pt
{\bf 3.1.} {\it Right--universal modules.}
Let $A$ be an algebra and $W$ a $\bk$-module. Consider
the $\bk$-module $W\ot A$. It has a canonical right $A$-module
structure $(w\ot x).y=w\ot xy$ and with this structure
$W$ generate $W\ot A$.
Observe that in the case $\bk$ is a field $W\ot A$ is
a free right module and $W$ plays the role of space spanned by its 
free generators. Universal bimodules of the form $A\ot W\ot A$
has been considered in [8].

Notice the following universal property of $W\ot A$ : 
for each $\bk$-linear mapping $\Phi$ from $W$ 
into a right $A$-module $M$ there exists a unique  $A$-right
module homomorphism $\Phi^u:W\ot A\ra M$ (called an {\it
universal lift} of $\Phi$), such that
$$\Phi^u\c i=\Phi  \leqno(9)$$
where, $i(w)=w\ot 1_A$ is a canonical mapping $i:W\ra W\ot A$.
$\Phi^u$ is defined by setting $\Phi^u(w\ot x)=\Phi(w).x$ and its 
image is the right submodule generated by $W$. Therefore, $\Phi^u$ is 
surjective if and only if $\Phi(W)$ generates $M$. In other words 
each right $A$-module generated by $W$ is a factor module
of $W\ot A$. More exactly we have

\th{Lemma}{3.1.1.}{Let  $W$ be a $\bk$-submodule of a right $A$-module $M$.
$W$ generates $M$ if and only if
$$M\cong (W\ot A)/\ker \id_W^{\,u}$$
where, $\id_W^{\,u}$  is an universal lift of the identity map
$\id_W:W\ra W\subseteq M$.
}
\vskip4pt plus2pt
If \ $\ker \id_W^{\,u}=0$ then $W$ is said to be a space of
{\it universal generators} of $M$. In this case $M\cong W\ot A$ and is
said to be a {\it right--universal} module generated by $W$.
$W$ is said to be a {\it $\bk$-module of universal generators} of $M$.
\remar{Note\ }{
We do not assume  that the $\bk$-module $W$ 
possesses free generators, {\it i.e.} that $W$ is a free
$\bk$-module. Therefore, $W\ot A$ is not a free right $A$-module,
in general. However, starting from the Subsection 3.4 we
shall assume projectivity of $W$.}
The map $u$ which establishes the correspondence between a $\bk$-linear
maps from $W$ into $M$ and an $A$-linear maps from $W\ot A$ into $M$ 
$$
\hom(W, M)\longleftrightarrow \hom_A(W\ot A, M) \leqno(10)
$$
is $\bk$-linear and bijective.
\vskip4pt plus2pt
{\bf 3.2.} {\it Left $A$--module structure.}
Since for calculi we need bimodule structure rather then right
module we investigate now how to determine a left module structure
on $W\ot A$ in such a way that $W\ot A$ becomes a bimodule.

Recall that a left $A$-module structure on a $\bk$-module $M$ is
given by an arbitrary element of $\alg (A, \nd M)$. Due to a canonical
isomorphism of $\bk$-modules (see [2] p. 267) 
$$\hom (A, \nd M)\cong\hom (A\ot M, M) \leqno(11)$$
one has $\alg (A, \nd M)\subset\hom (A\ot M, M)$.

When $M$ is a right $A$-module then we are interested in such left
module structures which convert $M$ into an $A$-bimodule. The set of
such structures is now identical with $\alg (A, \nd_A M)\subset
\hom_A (A\ot M, M)$.

Let $V,\ W$ be two $\bk$-modules. In this note by a {\it twisting map}
or simply a {\it twist} we shall call an arbitrary $\bk$-linear
map from $V\ot W$ into its twist $W\ot V$.

Proceeding to the case of a right-universal $A$-module $W\ot A$ one sees
by the universal mapping property (10) that
$\hom_A(A\ot W\ot A, W\ot A)$ is isomorphic to $\hom (A\ot W, W\ot A)$.
Consider arbitrary twisting map $b:A\ot W\ra W\ot A$. We will be 
interested to know  when its universal lift $b^u :A\ot W\ot A\ra W\ot A$ 
determines an $A$-bimodule structure on 
the right--universal $A$-module $W\ot A$.
Thus
\th{Lemma}{3.2.1.}{Let $A$ be an algebra and $W$ a $\bk$-module.
A twisting map $b:A\ot W\ra W\ot A$ defines a bimodule structure
on $W\ot A$ if and only if the universal lift $b^u$ is an
algebra homomorphism i.e. $b^u\in \alg (A,\> \nd_A(W\ot A))$
}
\Proof One should remember that in this case
$\alg (A,\> \nd_A(W\ot A)) \subset \hom_A(A\ot W\ot A, W\ot A)$
by a canonical isomorphism mentioned above.
\vskip4pt plus2pt
\th{Proposition}{3.2.2.}
{There is a one to one correspondence between  left $A$-module
structures transforming $W\ot A$ into bimodule and  twisting
maps $b :A\ot W\ra W\ot A$ satisfying:
\item{(i)} $b(1_A\ot w)=w\ot 1_A$
\item {(ii)} $b\c (m\ot \id_W)= (\id_W\ot m)\c (b\ot \idA)\c (\idA\ot b)$
on \  $A\ot A\ot W$.}
\Proof $W\ot A$ already has a right $A$-module structure.
$b$ determines a left multiplication by $x.(w\ot y) = b(x\ot w).y
\equiv b^u(x\ot w\ot y)$.
Bimodule axioms follow easily from the postulates.
\vskip4pt plus2pt
Twistings satisfying the conditions (i) and (ii) above will be called
{\it bimodule commutation rules}  and the corresponding $A$-bimodule
structure will be denoted by $W\bot A$.
\th{Lemma}{3.2.3}{Let $M$ be an $A$-bimodule. Assume that the $\bk$-submodule
$W$ universally generates $M$ as a right $A$-module. Then there exists
a unique bimodule commutation rule $b:A\ot W\ra W\ot A$  such that
$M\cong W\bot A$.}
\Proof  Let $\Phi :W\ot A\ra M$ denotes an isomorphism of right
modules. Define $b(x\ot w)=\Phi^{-1}(x.w)$. The proof is done.
\vskip4pt plus2pt
$W\bot A$ is said to be a {\it right--universal bimodule} generated
by $W$ with a left module structure defined by the bimodule commutation
rule $b$.

Our next task is to investigate calculi on factor modules of the type
$W\ot (A/I)\cong (W\ot A)/(W\ot I)$, where $I$ is an ideal of $A$.
\remar{Note\ }{To be precise one has to distinguish between $\bk$-module
$W\ot I$ and its image $\im (W\ot I)\subset W\ot A$ under a canonical
map $\id_W\ot i$ where, $i:I\ra A$ is an injection (see [2] p. 252). 
Hence, $W\ot (A/I)\cong (W\ot A)/\im (W\ot I)$ in general.
For  a projective $\bk$-module $W$ it is possible to identify 
$W\ot I$ with its image.
}

With this notation one has
\th{Proposition}{3.2.4. \ }{ Let $I\,<\, A$ be an ideal of $A$ and
$W$ be a projective $\bk$-module.
The bimodule commutation rule $b:A\ot W\ra W\ot A$ factorize to
the bimodule commutation rule  $\tilde b:(A/I)\ot W\ra W\ot (A/I)$
if and only if
$$
b(I\ot W)\subseteq W\ot I .\leqno (12a) $$
If in addition $d:A\ra W\bot A$ is a calculus then the condition
$$
dI\subseteq W\ot I \leqno(12b)$$
is sufficient and necessary for existence of a factor calculus
$\tilde d:A/I\ra W\,_{\tilde b}\!\ot (A/I)$.
}
\Proof First consider factor bimodule structure.
Take $L=W\ot I$. We have to prove that $L$ is an $A$- subbimodule
of $W\bot A$ (Lemma 2.1.1) which satisfies the condition (5) of 
Proposition 2.1.2 . It is clear that $L.A\subseteq L$ and $M.I\subseteq L$.
It means that the factor $\bk$-module $W\ot (A/I)\cong (W\ot A)/L$
has a right $A/I$-module structure 
and $W$ universally generates it as a right $A/I$ module.
Next calculate $A.L=b^u (A\ot W\ot I)=b(A\ot W).I\subset 
(W\ot A).I\subseteq L$. To see $I.M\subseteq L$ one has to use (12a).
Indeed $I.M=b^u(I\ot W\ot A)=b(I\ot W).A\subseteq (W\ot I).A\subseteq L$.
Therefore, we have a bimodule structure on a factor module $W\ot (A/I)$.
Then by Lemma 3.2.3 it is isomorphic to $W\,_{\tilde b}\!\ot (A/I)$
with some  commutation rule $\tilde b :(A/I)\ot W\ra W\ot (A/I)$.
In this case 
$$
\tilde b((x+I)\ot w)=b(x\ot w) + L.
$$

Now the last statement is a simple consequence of Theorem 2.2.1 .
\vskip4pt plus2pt
The condition (12a) will be called a  compatibility condition 
between an ideal $I$ and a bimodule commutation rule $b$. 
The conditions (12a) and (12b) form together 
consistency conditions for the existence of a factor calculus
with values in right--universal bimodules.

\remar{Remark\ {3.}\ }{Notice that in the conditions (12a), (12b) both 
sides depend on the ideal $I$. Therefore the existence of an optimal 
algebra is not in general an obvious question [3--6].}

There exists a canonical $\bk$-module homomorphism
$$
\lambda :\nd (W)\ot A\longrightarrow \hom (W,  \> W\ot A)
$$
which is in general not surjective nor injective ([2] p. 275 formulae (22)).
It induces a homomorphism $\rho=\hom (\id_A, \lambda)$ ([2] p. 229 and 267):
$$
\rho:\hom (A, \> \nd (W)\ot A)\longrightarrow \hom (A, \> \hom (W, \> W\ot A))
\cong\hom (A\ot W, \> W\ot A)
$$
assigning to each $\bk$-linear map $b:A\ra \nd (W)\ot A$
a twist $b^\rho :A\ot W\ra W\ot A$ by the formulae
$$b^\rho(x\ot w)=b(x)(w)\in W\ot A  \  
$$
where, $b^\rho\equiv \rho(b)$. It is not difficult to check that the 
image $\rho (\alg (A, \> \nd (W)\ot A))\subset \hom (A\ot W, \> W\ot A)$ 
provides bimodule commutation rules. Therefore, the compatibility
condition (11) can be replaced in this case by the following condition
$$
b^\rho (I)\subseteq \nd (W)\ot I \ . \leqno(13)$$

$\lambda$ is bijective (it is so if {\it e.g.} $W$ is finitely generated 
and projective $\bk$-module [2] p. 274) if and only if $\rho$ is bijective 
([2] p. 229). This approach is more convenient for the case when $\bk$ 
is a field (c.f. [3--6]).

{\bf 3.3.} {\it Calculi on tensor algebra.}
Let $V$ be a $\bk$-module and $TV$ denotes its tensor algebra.
(We neglect a graded algebra structure of $TV$ at this stage).
It is a free $\bk$-algebra generated by $V$ since
it possesses the following universal property: for each linear
map $g$ from $V$ into algebra $A$ there exists a unique algebra 
homomorphism $g_u: TV\ra A$ (called again an {\it universal lift} 
of $g$), such that 
$$g_u\c j=g $$
where, $j$ denotes a canonical
inclusion of $V$ into $TV$. The image of $g_u$ is a subalgebra generated
by $g(V)$. In particular, if $V\subseteq A$ generates $A$ then
$A\cong TV/I$ where, $I=\ker g_u$ is an ideal of relations in $A$.
In other words any algebra generated by $V$ is a quotient of $TV$.
Again the correspondence between $\bk$-linear maps from $V$ into
$A$ and algebra maps from $TV$ into $A$ {\it i.e.} between
$$\hom(V, A)\longleftrightarrow \alg (TV, A) \leqno(14)$$
is an isomorphism of $\bk$-modules.
\remar{Note\ }{Let $I<TV$ be a two-sided ideal. It appears that in general 
the factor algebra $TV/I$ does not contain $V$ as a $\bk$-submodule.
Some of the results below are valid in the general situation $TV/I$.
Some other, require an additional assumption that $V\subseteq TV/I$ .
In this case we shall speak that algebra is generated by its
submodule $V$. Here $V$ need not to be a free $\bk$-module.
Therefore, $TV$ is not a free algebra in general.
}
From now on, up to the end of this note we shall assume
that $V$ and $W$ are projective $\bk$-modules.
In a case of $V$-universal algebra $TV$ one has (c.f. [4]) 
\th{Proposition}{3.3.1}{Any (linear) map $b_0:V\ot W\ra W\ot TV$
determines uniquely a bimodule structure on $W\ot TV$.
Conversely, any bimodule commutation rule on $W\ot TV$ can be
obtained in this way.}
\Proof We shall make use of canonical isomorphisms
listed above. First by the canonical lift (10) we get an element
$(b_0)^u\in \hom_{TV}(V\ot W\ot TV,  \> W\ot TV)$ which is a $TV$-right
module map. Next, by a canonical isomorphism (11)
$(b_0)^u\in \hom (V, \> \nd_{TV}(W\ot TV))$. Finally, a canonical lift 
(14) gives an algebra map $((b_0)^u)_u\in \alg (TV, \> \nd_{TV}(W\ot TV))$.
This is a required left $TV$-module structure on $W\ot TV$. By Lemma 3.2.1
it leads to the bimodule commutation rule $b: TV\ot W\ra W\ot TV$
such that $b^u =((b_0)^u)_u$. It finishes the proof.
\vskip4pt plus2pt
\remar{Remark\ {4 .}\ }{Observe that since $TV=\oplus_{n\geq 0}V^{\ot n}$
is a graded algebra then $W\ot TV\equiv \oplus_{n\geq 0}W\ot V^{\ot n}$
is a graded right module. Therefore, $b$ preserve gradation imply
that $W\bot TV$ is a graded bimodule. If in addition $I$ is a graded
ideal then factor bimodule is again graded. Replace word 'graded'
by 'filtered'. The same is true.
}
\th{Proposition}{3.3.2.\ (Calculus on tensor algebra)}
{For each two $\bk$-linear maps $d_0:V\ra W\ot TV$ and
$b_0:V\ot W\ra W\ot TV$ there exists exactly one 
calculus $d$ on $TV$ with values in $W\bot TV$ such that
$d|_V=d_0$ and $b|_{V\ot W}=b_0$. Moreover in the graded case
i.e. $b_0:V\ot W\ra W\ot V$ and $d_0:V\ra W\ot 1_A$,
the derivation $d$ is a graded map of degree -1.}
\Proof (c.f. also [3, 4]). 
Let $M$ be a $TV$-bimodule and $d_0:TV\ra M$ a $\bk$-linear map.
One has to define an $M$-valued $\bk$-derivation $d$ on $TV$. $V$ generates 
$TV$. By  $\bk$-linearity it is enough to 
define $d$ on monomials $v_1\ot v_2\ot\ldots\ot v_k \in TV$ where,
$v_1, v_2, \ldots , v_k\in V$. Do it by induction.
\item{(i)} $dv=d_0(v)$ for $v\in V$.
\item{(ii)} If $x=v\ot x_1$  is a word of length $k$, $v\in V$ and
$x_1$ is a word of length $(k-1)$ then define 
$$
dx=  (dv). x_1 + v. dx_1 = (d_0 v). x_1 + v. dx_1).
$$

Now, we have to prove that $d$ satisfies the Lebniz rule.
Let $x, y$ be two monomials and $x=v\ot x_1$ as before. Again by
inductive assumption if $d(x_1\ot y)=(dx_1).y + x_1.dy$ then
$$
d(x\ot y)=d(v\ot (x_1\ot y))= (d_0v). (x_1\ot y) + v. d(x_1\ot y) =
$$
$$
=((d_0v). x_1). y + v. (dx_1. y + x_1. dy) = dx. y +  x. dy
$$
since $v. (x_1. dy) = (v\ot x_1). dy = x. dy$.
Thus we have prooved that any $\bk$-linear map $d_0: V\ra M$ uniquely extends
to a $\bk$-derivation $d: TV\ra M$.
Obviously, any calculus on $TV$ can be obtained in this way.

Next, apply above result to the case $M=W\bot RV$ where,
$b: TV\ot W\ra W\ot TV$ denotes the
commutation rule constructed out of $b_0$ in the proof of Proposition 3.3.1 .

Finally, in the homogeneous case, $d_0v\in \bk$ and $b_0:V\ot W\ra W\ot V$.
In this case $b$ preserve the grading (which is a length of monomials).
Therefore, by inductive formulae $d$ decreases a degree by one.
The proof is done.
\vskip4pt plus2pt

In other words calculus on tensor algebra with values in right universal
bimodule is completely determined by the {\it initial data}: $b_0$ and $d_0$.
\th{Proposition}{3.3.3.\ }{
Let calculus $d:TV\ra W\bot TV$
is given. The condition  
$$
d I + b (I\ot W)\subseteq W\ot I \ . \leqno (15)$$
is a sufficient and necessary
for the existence (a unique) factor calculus 
$\tilde d:TV/I\ra W\,_{\tilde b}\!\ot (TV/I)$
where, $\tilde b$ denotes a factor commutation rule.
}
\Proof It is a simple consequence of Proposition 3.2.4 .
\vskip4pt plus2pt
\remar{Remark\ {5.}\ }{It is clear that the condition (15) 
can be equivalently replaced by two conditions
$$d I \subseteq W\ot I \leqno(15a)$$
$$b (I\ot W)\subseteq W\ot I \ \ (\ \Longleftrightarrow
b ^\rho (I)\subseteq \nd (W)\ot I\ ) \leqno(15b)$$
(Formulae in brackets relates to the case governed by (13)).
For the quadratic algebras this type of conditions, when written in the 
operator form (see below) have been introduced in [17]. 
Therefore, they are often called Wess--Zumino linear and quadratic 
consistency conditions. In some cases it is possible to express the
consistency conditions in terms of generating $\bk$ -module for the
ideal $I$. It will be studied in Subsection 3.5 .
}
One sees from the above considerations that factor calculi 
of calculi on universal algebra form
an important class of calculi. This is so due to the following
two facts. Firstly, we are able to describe all possible calculi
on universal algebra with values on right--universal bimodule. 
Secondly, the consistency conditions (15a, b) are convenient and more 
easy to handle. In the next point we answer the question:
under which circumstances all calculi on the algebra $A$
can be obtained in this way.
\vskip4pt plus2pt
{\bf 3.4.} {\it Calculi with right--universal bimodule -- structure
theorem.} Recall that $V$ is a projective $\bk$-module if and only if for
every surjective ($\bk$-linear) map $\Pi :P\ra Q$ and every  map
$f:V\ra Q$ there exists a map $\hat f:V\ra P$ such that 
$\hat f\c\Pi=f$.
It is said that $f$ can be lifted to a map of $V$ into $P$ or
that $\hat f$ covers $f$ subject to the epimorphism $\Pi$.
Of course, $\hat f$ exists but is not unique in general. 
Tensor product of projective modules (over commutative ring) 
is again projective.

Joining the universal properties of $W\ot A$ with that of $TV$ one has:
Let $A$ be an algebra generated by $\bk$-submodule $V\subseteq A$ and
$M$ be a right--universal $A$-bimodule universally generated by 
its submodule $W\subseteq M$. Then $A\cong TV/I$ and 
$M\cong W\bot (TV/I)$ for some bimodule commutation rule $b$.
We show our main result that  each calculus on $TV/I$ with values in 
right--universal bimodule $W\bot A$ is a factor calculus of some
{\it cover} calculus $\hat d$ on tensor algebra $TV$ subject to the cover
bimodule commutation rule $\hat b: TV\ot W\ra W\ot TV$. More exactly one has

\th{Theorem}{3.4.1.\ (Main result. Cover calculus.)}{
With the notation as above :

Let $d: TV/I\ra W\bot (TV/I)$ be a differential 
calculus. There exist (not unique in general)  initial data 
$\hat{b}_0:V\ot W\ra W\ot TV$ and $\hat{d}_0:V\ra W\ot TV$ 
such that the corresponding
calculus $\hat d$ on $TV$ factorize to $d$. Thus the diagram
$$
\matrix{TV &{\buildrel \hat{d}\over\longrightarrow} & W\,_{\hat b}\!\ot TV\cr
\pi \downarrow &   & \downarrow \id_W\ot\pi\cr
TV/I &{\buildrel d\over\longrightarrow} & W\bot (TV/I) 
}$$
is commutative where, $\pi :TV\ra TV/I$ is a canonical surjection.
Moreover, the following consistency condition is satisfied
$$
\hat d I + \hat b (I\ot W)\subseteq W\ot I \ . \leqno (16)$$
}
\Proof Since $\pi$ and $\id_W\ot\pi$ are in particular surjective
$\bk$-linear maps one can apply projectivity of $V$ and $V\ot W$. Set  
$b_0=b|_{V\ot W}:V\ot W\ra W\ot A$ and $d_0=d|_V:V\ra W\ot A$.
Let $\widehat{b_0}:V\ot W\ra W\ot TV$ and $\widehat{d_0} :V\ra W\ot TV$
be some their (not unique) lifts. By Proposition 3.3.2 let $\hat d$ denotes 
the corresponding (unique) calculus on $TV$. To show commutativity of
the diagram ( it is enough to check it on generating space $V$
(Lemma 1.2.2). We leave it to the reader (c.f. [3, 4]).
The consistency condition (16) follows easily from Proposition 3.3.3 .
\vskip4pt plus2pt
\remar{Remark\ {6.}\ }{So called coordinate calculus [3--6]
corresponds to the case $W=V$ be a vector space and $d_0=\id_V$. 
}
Theorem shows the way to construct and classify all
calculi with right--universal bimodules if one has been given the 
$\bk$-projective modules of generators of an algebra and a bimodule. 
In order to do this we have to
follow one of two main strategies:
\item{-} fix the ideal $I$ in the tensor algebra and therefore an algebra
$A$ and its presentation $A\equiv TV/I$. Look for all 
initial data $d_0,\ b_0$ such that the corresponding calculus $d$ on
$TV$ satisfies the consistency conditions (15a, b);
\item{-} choose the initial data. Look for all consistent ideals $I$.
\vskip4pt plus2pt
{\bf 3.5.}  {\it Applications. Calculi on quadratic algebras.}
Let $V$ be a $\bk$-module and $K\subseteq V^{\ot 2}$ be a submodule.
Denote by $<K>$ an ideal generated in $TV$ by $K$. The factor
algebra $TV/<K>$ is said to be an {\it algebra with quadratic relations}
({\it quadratic algebra} in short). It is not difficult to see that the
ideal $<K>$ is a graded ideal $<K>=\oplus_{n\geq 2}K_n$ with the
homogeneous components satisfying the recurrent relations
$$
K_{n+1} = V\ot K_n + K_n\ot V . \leqno(17) $$
where, $K_2\equiv K$. Let $W$ be another $\bk$-module. 
Assume we are given the derivation $d:TV\ra W\ot TV$ with
homogeneous initial data $b_0:V\ot W\ra W\ot V$ and $d_0:V\ra W$.
In this case the derivation $d$ decrease the degree. Therefore,
the linear consistency condition (15a) takes the very simple form
$$
d K=0 \ .\leqno(18)$$
Taking into account the Leibniz rule and homogeneity of $d$ we get
$$[d_0\ot\idV + b_0\c (\idV\ot d_0)]|_K = 0 \ .\leqno(19)$$
Also quadratic consistency condition (15b) simplify to
$$
b(K\ot W)\subseteq W\ot K \ (\ \Longleftrightarrow \ 
b^\rho (K)\subseteq \nd(W)\ot K\ ) \leqno(20)$$
where, $b:TV\ot W\ra W\ot TV$ is an associated bimodule commutation rule.
To see this one has to apply algebra homomorphism $b^\rho$ to 
formulae (17).
It appears  that (20) can be further simplify to
$$
(b_0\ot\idV)\c (\idV\ot b_0) (K\ot W)\subseteq W\ot K  \leqno(21)$$
if we use Proposition 3.2.1 . 

Our aim in this Section is to find out an 
operator (matrix) form of consistency condition (19, 21). Thus one has
(c.f. [5, 12])
\th{Proposition}{3.5.1.}{Let $c\in\nd (V^{\ot 2})$ be a twist with the 
property that $\im\, c\equiv c(V\ot~V)=K$.
Then the consistency condition (19, 21) for the homogeneous calculus on the
quadratic algebra $TV/<\im\, c>$ are equivalent to the following conditions
$$[d_0\ot\idV + b_0\c (\idV\ot d_0)]\c c = 0 \leqno(22)$$
$$(b_0\ot\idV)\c (\idV\ot b_0)\c (c\ot\id_W)=(\id_W\ot c)\c a\leqno(23)$$
for some twist $a:V^{\ot 2}\ot W\ra W\ot V^{\ot 2}$.}
\remar{Example\ {}}{
Let $b_0:V\ot W\ra W\ot V$ be simply a switch {\it i.e.} $b_0(v\ot w)=w\ot v$.
Then for any ideal $I < TV$ the quadratic consistency condition (15b) is 
automatically satisfied. It follows that for any $\bk$-algebra $A\equiv TV/ I$
generated by $V$, one has $b(x\ot w)=w\ot x$ 
and $x.(w\ot y).z=w\ot (xyz)$ in $W\bot A$. The linear consistency
condition (15a) is the only one to be fulfilled for the existence of factor
calculus on $A$. In the quadratic algebra case with $d_0:V\ra W$ 
(Proposition 3.5.1)\ \ (22) reads
$$
(d_0\ot\idV)\c c\ =\ 0 \leqno(24)$$
{\it i.e.} $\im c\subset \ker (d_0\ot\idV)$. Therefore, $TV/\,<\ker 
(d_0\ot\idV)>$ is an optimal algebra for these bimodule commutation rules.}
\remar{Remark\ {7.}\ }{
Very peculiar solutions of the quadratic consistency condition (23) are
given by solutions of the braid equation (quantum Yang--Baxter eq. ) vis.
$$
(b_0\ot\idV)\c (\idV\ot b_0)\c(c\ot\id_W)=(\id_W\ot c)\c 
(b_0\ot\idV)\c (\idV\ot b_0) \leqno(25)
$$
on $V^{\ot 2}\ot W$. (See [5] for more examples and comments.)}
\vskip4pt plus2pt
\references{ref}{
\item{[1]} H. C. \spa{Baehr}, A. \spa{Dimakis} and F. \spa{
M{\"u}ller-Hoissen}, {\it Differential Calculi on Commutative Algebras}\/,
J. Phys. A: Math. Gen. 28 (1995), 3197--3222, hep-th/9412069.
\item{[2]} N. \spa{Bourbaki},
{\it Elements of mathematics.Algebra I.Chapters 1--3},
Springer--Verlag, Berlin, 1989.
\item{[3]} A. \spa{Borowiec}, V. K. \spa{Kharchenko} and Z. \spa{Oziewicz},
{\it On free differentials on associative algebra}\/,
in~: Non-Associative Algebra and Its Applications, S. Gonz\'alez (ed.),
Kluwer Academic Publishers, Dordrecht 1994, ISBN
0-7923-3117-6, Mathematics and its Applications, vol. 303, 46--53,
(hep-th/9312023).
\item{[4]} A. \spa{Borowiec} and V. K. \spa{Kharchenko},
{\it Algebraic approach to calculi with partial derivatives}\/,
Siberian Advances in Mathematics 5, 2 (1995), 10--37.
\item{[5]} A. \spa{Borowiec} and V. K. \spa{Kharchenko},
{\it Coordinate calculi on associative algebras\rm,
in~: Quantum Group, Formalism and Applications}\/, J. Lukierski, 
Z. Popowicz and J. Sobczyk (ed.), Polish Sci. Publ. PWN Ltd., 
Warszawa, 1995. ISBN 83-01-11770-2, 231--241, q-alg/9501018.
\item{[6]} A. \spa{Borowiec} and V. K. \spa{Kharchenko},
{\it First order optimum calculi}\/, 
Bull. Soc. Sci. Lett. {\L}\'od\'z v. 45,
Ser. Recher. Deform. XIX, (1995), 75-88, q-alg/9501024.
\item{[7]} A. \spa{Borowiec}, {\it Cartan Pairs}\/,
Czech. J. Phys.  -- to be published, q-alg/9609\ \ .
\item{[8]} J. \spa{Cuntz} and D. \spa{Quillen}\/,
{\it Algebra Extension and Nonsingularity}, J. Amer. Math. Soc.
8, 2 (1995), p. 251--289.
\item{[9]} A. \spa{Dimakis}, F. \spa{M{\"u}ller-Hoissen} and
T. \spa{Striker},
{\it Non-commutative differential calculus and lattice gauge theory}\/,
J. Phys. A: Math. Gen. 26 (1993), 1927--1949.
\item{[10]} G. \spa{Maltsiniotis}, {\it Le Langage des Espaces et des 
Groupes Quantiques}\/, Commun. Math. Phys. 151 (1993), 275--302.
\item{[11]} Yu. I. \spa{Manin},
{\it Notes on quantum groups  and quantum de Rham complexes.}\/,
Preprint, MPI/91-60 (1991). 
\item{[12]} E. E. \spa{Mukhin},
{\it Yang--Baxter operators and noncommutative de Rham complexes}\/,
Russian Acad. Sci. Izv. Math. 58, 2 (1994), 108--131 (in Russian).
\item{[13]} R. S. \spa{Pierce}, {\it Associative algebras}\/, Graduate
Texts in Mathematics \# 88, Springer-Verlag, New York, 1982.
\item{[14]}  W. \spa{Pusz} and S. \spa{Woronowicz},
{\it Twisted second quantization}\/,
Reports on Mathematical Physics 27, 2 (1989), 231--257.
\item{[15]}  W. \spa{Pusz},
{\it Twisted canonical anticommutation relations}\/,
Reports on Mathematical Phys. 27, 3 (1989), 349--360.
\item{[16]}   K. \spa{Schm{\"u}dgen} and A. \spa{Sch{\"u}ler},
{\it Classification of bicovariant calculi on quantum spaces and
quantum groups}\/, C. R. Acad. Sci. Paris 316 (1993), 1155--1160.
\item{[17]} J. \spa{Wess} and B. \spa{Zumino},
{\it Covariant differential calculus on the quantum hyperplane}\/,
Nuclear Physics 18 B (1990), 303--312, Proc. Suppl. Volume in honor 
of R. Stora.
\item{[18]} S. L. \spa{Woronowicz},
{\it Differential calculus on compact matrix pseudogroups 
(quantum  groups)}\/, Commun. Math. Phys. 122 (1989), 125--170.
}

\bye